\documentclass[11pt]{article}
\usepackage[centertags]{amsmath}
\usepackage{amsfonts}
\usepackage{newlfont}
\usepackage{graphicx}
\usepackage{eurosym}
\usepackage[pdfpagelabels]{hyperref}
\usepackage{epstopdf}
\usepackage{amsmath,amssymb, xcolor}
\usepackage{caption}
\usepackage{subcaption}
\usepackage{amsfonts}
\usepackage{textcomp}
\usepackage{appendix}
\usepackage{newlfont}
\usepackage{graphicx}
\usepackage{eurosym}
\usepackage{bm}
\usepackage{mathtools}

\newcommand{\E}{\mathbf{E}}

\def\cF{{\cal F}}

\begin{document}

\title{{A note on D-functions and P-covariances on Hilbert spaces and related inequalities}}

\author{S. Scarlatti\thanks{Dept. of Economics and Finance - University of Rome Tor Vergata, \texttt{sergio.scarlatti@uniroma2.it}}}

\maketitle

\date{}
\begin{abstract}
In this note we first review the concept of D-function, closely connected with Cauchy-Schwarz inequality,  and then introduce the notion of P-covariance on a Hilbert space, where $P$ is an orthogonal projection. We show that when P is specialized to be one-dimensional many well-known inequalities such as Buzano, Richard  and Walker inequalities  are simple consequences of P-covariance inequalities and their relation with D-functions. By means  of these concepts  enhancements of the previous mentioned inequalities are also established with minimum effort. A more thorough analysis of Walker inequality is presented jointly with some novel financial considerations as well as a new refinement of Holder inequality. \medskip

\end{abstract}

\maketitle
\section{Introduction}
In Mathematics, so as in all applied sciences, quantitative comparisons of similar or even different quantities play a crucial role. Mathematical inequalities are key analytical results which express such comparisons in a neat way. In this note we focus on some classical inequalities in Hilbert spaces, but in the last part of it we shall also consider inequalities in $L^p$ spaces ($p>1$). 
We highlight some simple mathematical structures which underly a bunch of inequalities , some of which are very famous: Cauchy-Schwarz [1], Buzano [2], Richard [3], Walker [4]. The main aim of the work is to show how these structures, namely D-functions and P-covariances, allow for a fast-to-prove and  easy-to-generalize  type of approach to such inequalities. We also discuss some previously unnoticed aspects of Walker inequality related to the theory of financial investments [5], and use Walker to improve over classical Holder estimate on function spaces. 
\section{D-functions}
 Let $(H, \langle \cdot,\cdot\rangle_H )$ a real Hilbert space, $V$ any closed subspace of $H$  and $P\equiv P_V$ the corresponding orthogonal projection on $V$,  that is $PH=V$. For  any $x\in H$ we introduce the quantities $p(x)\equiv ||Px||$, $q(x)\equiv ||P^{\perp}x||$ and consider the two-dimensional positive vector
$$
v_P(x)\equiv (p(x),q(x));
$$
the following definition was set in [6]:\\
{\bf{Definition 1}} For fixed $P=P^2$ define on $H\times H$ the positive symmetric function
\begin{equation}\label{D-function}
D(x,y|P)\equiv \langle v_P(x),v_P(y)\rangle_{\mathbb{R}^2},
\end{equation}
that is $D(x,y|P)=||Px||||Py||+||P^{\perp}x||||P^{\perp}y||$, with $P+P^{\perp}=Id$.\\ We shall call such a function a $P$-decoupling function,  more briefly a decoupling function  [6], $P$ being implied , since it completely separates operations on $x$ from those on $y$.\\ 
A decoupling function, or D-function, has the following interesting properties:\\
{\bf{Proposition 1.}}\\
\noindent
(i)$D(x,x|P)=||x||^2$ \text{and} $D(x,y|P^{\perp})=D(x,y|P)\geq 0$,\\
(ii)$D(\lambda x,\mu y|P)=|\lambda\mu| D(x,y|P)$\; $\forall \lambda , \mu \in \mathbb{R}$ \;\text{and} $D(x+x',y|P)\leq D(x,y|P)+ D(x',y|P)$,\\
(iii)$D(P^\#x,y|P)=||P^\#x||\;||P^\#y||$ for $P^\#\in\{P,P^{\perp}\}$,\\
(iv)$D((P-P^{\perp})x,y|P)=D(x,y|P)$,\\
(v)$D(x,y|P)= 0 \;\text{for all} \;(x,y)\in V\times V^{\perp}$, with $V=PH$,\\
(vi) $D(x,y|P)\geq |\langle x,y \rangle_H|$\; ($D$-inequality),\;\\
(vii)for all $x,y \in H$ it holds
\begin{equation}\label{newCS}
D(x,x|P)D(y,y|P)-D^2(x,y|P)=det^2
\begin{pmatrix}
v_P(x)\\
v_P(y)
\end{pmatrix};
\end{equation}
notice that for $y=\alpha x$ the r.h.s. of (\ref{newCS}) has zero value, however this condition is only a sufficient one for having zero.\\
From the first and the last two properties it follows that for any $P=P^2$ we have
\begin{equation}\label{improved-CS}
|\langle x,y \rangle_H|\leq D(x,y|P)\leq ||x|| ||y||\;\forall x,y\in H 
\end{equation}
enhancing classical Cauchy-Schwarz inequality (CS). In particular the upper bound follows from (\ref{newCS}), or even more directly, from  classical CS-inequality applied to $ \langle v_P(x),v_P(y)\rangle_{\mathbb{R}^2}$. The lower bound is simply property (vi) which is mainly a consequence of orthogonal decomposition in $H$, that is $Id=P+P^{\perp}$ with $PP^{\perp}=0$, and of the classical CS-inequality [6]. This last fact is often reported as the "self-improvement" property.  We notice that considerations similar to those contained in [6] already appeared in [7], although the concept of D-function was not present there.\\
{\bf{Remark 1:}} We see that properties (vi) and (vii) automatically imply that for any projection $P$ the following inequality holds:
\begin{equation}\label{Dragomir}
||x||^2||y||^2-|\langle x,y \rangle|^2_H\geq det^2
\begin{pmatrix}
v_P(x)\\
v_P(y)
\end{pmatrix},
\end{equation}
generalizing a result proven in [8, (Thm.2)]  for one dimensional projections.\\
\section{P-covariances, Richard and Buzano}
{\bf{Definition 2.}}: Let $P=P^2$ be an orthogonal projection, we define on $H\times H$ the following symmetric bilinear form 
\begin{equation}\label{P-cov}
cov^P_H(x,y)\equiv \langle x-Px,y-Py\rangle_H=\langle x,y\rangle_H-\langle Px,y\rangle_H,
\end{equation}
for each pair $(x,y)$ we call the value of the form the P-covariance  between $x$ and $y$. We also name  the quantity $var_H^P(x)\equiv cov^P_H(x,x)=||x||^2-p(x)^2=q(x)^2\geq 0$ the P-variance of $x$. \\
For a fixed $z\in H$ of unitary norm let $P_z$ denotes the one dimensional orthogonal projection over the subspace $V\equiv\{\lambda z: \lambda\in \mathbb{R}\}$, with  $P_zx=\langle x,z\rangle_Hz$. For this special case we set the notation
\begin{equation}\label{z-cov}
cov^z_H(x,y)\equiv cov^{P_z}_H(x,y)=\langle x,y\rangle_H-\langle x,z\rangle_H\langle y,z\rangle_H;
\end{equation}
and call this number the $z$-covariance between $x$ and $y$. Clearly the above definition is such that it reproduces  classical covariances when specialized to function spaces, see e.g.  [10], or to random variables over probability spaces, see next section. However this generalization turns out to be useful, indeed we have:\\
{\bf{Proposition 2.}} The following inequalities hold for any projection $P$ and any pair $(x,y)$:\\
$$
|2\langle Px,y\rangle_H-\langle x,y\rangle_H|\leq D(x,y|P),\;\;\;\;\;\;\;\;(eR)
$$
$$
2|\langle Px,y\rangle_H| \leq D(x,y|P)+| \langle x,y\rangle_H|, \;\;\;\;\;\;\;\;(eB)
$$
\begin{equation}\label{ineq}
| \langle x,y\rangle_H|\leq  ||P^{\perp}x||||P^{\perp}y||+|\langle Px,y\rangle_H|,\;\;\;(eD)
\end{equation}
where $(eD)$ stands for "enhanced D-inequality", the other capital letters referring to names.\\
{\bf{proof:}} P-covariances enjoy the following upper bound:
\begin{equation}\label{rho1}
|cov^P_H(x,y)|\leq q(x)q(y),
\end{equation}
Indeed
\begin{equation}\label{rho2}
|cov^P_H(x,y)|=| \langle x-Px,y-Py\rangle_H|\leq  ||P^{\perp}x||||P^{\perp}y||
\end{equation}
by classical CS inequality. By adding $p(x)p(y)=||Px||\;||Py||$ to both sides of (\ref{rho2}) and noticing 
$$
|\langle Px,y\rangle_H+\bigl(\langle Px,y\rangle_H-\langle x,y\rangle_H\bigr)|\leq ||Px||\;||Py||+|-cov^P_H(x,y)|
$$
 we obtain $(eR)$. Furthermore, being
\begin{equation}\label{chain}
\bigl|| \langle x,y\rangle_H|-|\langle Px,y\rangle_H|\bigr|\leq |\langle x,y\rangle_H-\langle Px,y\rangle_H|\leq ||P^{\perp}x||||P^{\perp}y||,
\end{equation}
the following two inequalities hold true:
$$
| \langle x,y\rangle_H|-|\langle Px,y\rangle_H|\leq ||P^{\perp}x||||P^{\perp}y||\;\;\;\;\;(I)
$$
$$
|\langle Px,y\rangle_H|-| \langle x,y\rangle_H|\leq ||P^{\perp}x||||P^{\perp}y||\;\;\;\;\;(II).
$$
We rewrite inequality (II) as
\begin{equation}\label{B1}
|\langle Px,y\rangle| \leq ||P^{\perp}x||||P^{\perp}y||+| \langle x,y\rangle_H|,
\end{equation}
from which we easily obtain $(eB)$. Moreover, we also rearrange inequality $(I)$ as 
$$
| \langle x,y\rangle_H|\leq  ||P^{\perp}x||||P^{\perp}y||+|\langle Px,y\rangle_H|,
$$
proving $(eD)$. Being $|\langle Px,y\rangle_H|\leq ||Px||\;||Py||$ we notice this last estimate is even tighter than the $D$-inequality $| \langle x,y\rangle_H|\leq D(x,y|P)\leq ||x|| ||y||$, however it is not of decoupling type. 
$\square$\\
{\bf{Remark 2:}} Let $||z||=1$, then inequality $(eB)$  improves and extends classical Buzano estimate 
$$
2|\langle x,z\rangle_H\langle y,z\rangle_H|\leq ||x||||y||+| \langle x,y\rangle_H|,\;\;\;\;(B)
$$
and inequality $(eR)$ improves and extends classical Richard inequality
$$
|2\langle x,z\rangle_H \langle y,z\rangle_H-\langle x,y\rangle_H| \leq ||x|| ||y||, \;\;\;(R).
$$
{\bf{Remark 3:}} In [9] it was proven that for any projection $P$ it holds $2|\langle Px,y\rangle_H| \leq ||x|| ||y|| +| \langle x,y\rangle_H|$, so $(eB)$ expresses a tighter bound.\\

It is remarkable how  easily and quickly such generalizations have been obtained by means of D-functions and $P$-covariances.\\

\section{ Walker, Sharpe and Holder}
By applying $(eD)$ to the case  $P\equiv P_z$, being $|\langle P_zx,P_zy\rangle_H|=||P_zx||\;P_zy||$, we just have
\begin{equation}\label{W}
| \langle x,y\rangle_H|\leq  ||P^{\perp}_zx||||P^{\perp}_zy||+|\langle x,z\rangle_H\langle y,z\rangle_H|=D(x,y|P_z).
\end{equation}
Walker inequality nicely fits the $z$-covariance framework corresponding to the choice of the following triplet: $H=L^2\equiv L^2(S,\cF,Q)$ ( $\cF$ a sigma-algebra, $Q$ a probability measure), inner product $\langle X,Y\rangle_{L^2} \equiv \E (XY)$ and  $z^*\equiv 1\equiv 1_{S}$ ( the characteristic function of a set) as fixed unitary vector.  By consequence
$$
cov_{L^2}^1(X,Y)=\langle X,Y\rangle_{L^2}-\langle X,1\rangle_{L^2}\langle Y,1\rangle_{L^2}=cov(X,Y),
$$
where $cov(X,Y)$ denotes classical statistical covariance (see also [6],[7]). We also set $\sigma^2_X\equiv var_{L^2}^1(X)=\E(X^2)-\E(X)^2$. For any $X\in L^2$ define the one-dimensional projection $P_1X\equiv \langle X,1\rangle_{L^2} 1=(\E X)1$ and the vector $v_P(X)\equiv (|\E(X)|,\sigma_X)$. We have:
\begin{equation}\label{seeH}
(\E(XY))^2\leq ||X||_2^2||Y||_2^2-\biggr(|\E X|\sigma_Y-|\E Y|\sigma_X\biggl)^2\leq ||X||_2^2||Y||_2^2 \;\;\forall X,Y \in L^2,
\end{equation}
however since, by (vii) of sect. 1, it holds
$$
D^2(X,Y|P_1)=||X||_2^2||Y||_2^2-\biggr(|\E X|\sigma_Y-|\E Y|\sigma_X\biggl)^2
$$
we have that (\ref{seeH}) is just a particular instance of (\ref{W}),  which in turn is an instance of $(eD)$. It is worth noticing that if $X$ and $Y$ are such that $|\E X|\sigma_Y=|\E Y|\sigma_X$, or equivalently, if it holds the relation
\begin{equation}\label{ratio}
\biggl(\frac{\E X}{\sigma_X}\biggr)^2=\biggr(\frac{\E Y}{\sigma_Y}\biggr)^2,
\end{equation}
the improvement over CS inequality vanishes. 
It is attractive to interpret condition (\ref{ratio}) also in financial terms. Indeed suppose that $X$ and $Y$ are future excess-returns of two different portfolios investing in the stock market, where excess-return means the extra return w.r.t the risk-free rate. Adopting a Markowitz's  investment style the numbers $\sigma_X$ and $\sigma_Y$ are interpreted as the respective risks of these two financial positions. Then the Sharpe Ratio (SR) of a financial position offering $X$ as future random return, introduced by W.Sharpe (who shared a Nobel price with H.Markowitz and M.Miller), is exactly given as the ratio $SR(X)\equiv \frac{\E X}{\sigma_X}$, a number which eventually can be negative if the investment performs badly. In general,  the higher the ratio the more attractive the investment. On the basis of (\ref{ratio}) we conclude that Walker inequality produce no improvement w.r.t. CS inequality if and only if $X$ and $Y$ are such that the squares of their Sharpe Ratios equalize. For more on this approach to investments and related questions see [5] and the paper [11].
\\


We now show that estimate (\ref{seeH}) can be used to improve Holder estimate:\\
 {\bf{Proposition 3.}} For $X\in L^p$ and $Y\in L^q$, $p,q>0$, $\frac1{p}+\frac1{q}=1$, it holds
 \begin{equation}\label{new-H}
\E(|XY|) \leq ||X||_p||Y||_q\biggl(\frac{1}{p^{2}}+\frac{1}{q^{2}}+\frac{2}{p q}\sqrt{1-\biggl(\frac{\sigma_V\E|X|^{p / 2}-\sigma_U\E|Y|^{q / 2}}{||X||_p^{p/2}||Y||_q^{q/2}}\biggr)^2}\;\;\biggr)\leq ||X||_p||Y||_q.
\end{equation}
{\bf{proof :}} by Young inequality it is easy to show that for $X\in L^p$ and $Y\in L^q$, $p,q>0$, it holds (see e.g. [12]):
$$
\E(|XY|) \leq\left(\frac{1}{p^{2}}+\frac{1}{q^{2}}\right)\|X\|_{p}\|Y\|_{q}+\frac{2}{p q}\|X\|_{p}^{1-\frac{p}{2}}\|Y\|_{q}^{1-\frac{q}{2}} \E(|X|^{p / 2}|Y|^{q / 2}) .
$$
Now we set 
$U\equiv |X|^{p / 2}\in L^2$ and $V\equiv |Y|^{p / 2}\in L^2$ and therefore, by (\ref{seeH}), we have
$$
\E(UV)\leq \sqrt{||U||_2^2||V||_2^2-\biggr(|\E U|\sigma_V-|\E V|\sigma_U\biggl)^2}=\sqrt{||X||_p^p||Y||_q^q-\biggr(\sigma_V\E|X|^{p / 2}-\sigma_U\E|Y|^{q / 2}\biggl)^2}
$$
$$
=||X||_p^{p/2}||Y||_q^{q/2}\sqrt{1-\biggl(\frac{\sigma_V\E|X|^{p / 2}-\sigma_U\E|Y|^{q / 2}}{||X||_p^{p/2}||Y||_q^{q/2}}\biggr)^2}
$$
Henceforth, for $\frac1{p}+\frac1{q}=1$, we have (\ref{new-H}).$\square$

Let us remark that this estimate differs  from other improvements appearing  in literature, see e.g. [13], at least for one aspect: in (\ref{new-H}) informations about $X$ and $Y$ enter in a completely separated way which allows the bound to be quickly evaluated using a parallel mode computation for the moments of $X$ and $Y$. Furthermore, in the case $p=q=2$ inequality (\ref{new-H}) reduces to
$$
\E(|XY|) \leq  ||X||_2||Y||_2\biggl(\frac1{2}+\frac1{2}\sqrt{1-\biggl(\frac{\sigma_Y\E|X|-\sigma_X\E|Y|}{||X||_2||Y||_2}\biggr)^2}\;\;\biggr)
$$
\begin{equation}\label{new-Walker}
=\frac{1}{2} ||X||_2||Y||_2+\frac1{2}\sqrt{||X||_2^2||Y||_2^2-\biggr(|\E X|\sigma_Y-|\E Y|\sigma_X\biggl)^2};
\end{equation}
Notice that this last estimate does not reproduce (\ref{seeH}) exactly: indeed in (\ref{seeH}) the first term of  (\ref{new-Walker}), that is the quantity $\frac1{2}||X||_2||Y||_2$, appears to be already replaced by the smaller quantity
$$
\frac1{2}\sqrt{||X||_2^2||Y||_2^2-\biggr(|\E X|\sigma_Y-|\E Y|\sigma_X\biggl)^2};
$$  
so we may say that (\ref{new-Walker}) contains just half of the improvement  given by (\ref{seeH}).
\section{Conclusions} This paper highlights the enhancements of some well-known and useful  inequalities which can be obtained with great simplicity and clarity by means of the concepts of D-functions and P-covariances over a Hilbert space.
In addition,  considerations of financial nature widen the interpretation of some technical condition entering in Walker inequality. The same inequality is then used to produce a new refinement of Holder inequality. Overall , the obtained results pave the way to further investigations on the possible range of applicability to information theory , computer science and physics of the presented improvements.\\
\noindent






\end{document}